\documentclass[12pt,twoside]{amsart}
\usepackage{amsfonts}
\usepackage{amsfonts}
\usepackage{amsfonts}
\usepackage{amsfonts}
\usepackage{amsfonts}
\usepackage{amsfonts}
\usepackage{amsfonts}
\usepackage{amsfonts}
\usepackage{amsfonts}
\usepackage{amsfonts}
\usepackage{amsfonts}
\usepackage{amsfonts}
\usepackage{mathrsfs}
\usepackage{mathrsfs}
\usepackage{txfonts}
\usepackage{bbm}
\usepackage{mathrsfs}
\usepackage{amsfonts}
\usepackage{amsmath}
\usepackage{amssymb}
\usepackage{wasysym}
\usepackage{psfrag}
\usepackage{graphics,epsfig,amsmath}  
\usepackage{comment}
\usepackage{geometry}
\usepackage{anysize}
\usepackage{fancyhdr}
\usepackage{indentfirst}
\usepackage{graphicx}
\marginsize{4cm}{4cm}{3.5cm}{3.5cm}

\newtheorem{Theorem}{Theorem}[section]
\newtheorem{Definition}[Theorem]{Definition}
\newtheorem{Lemma}[Theorem]{Lemma}

\newtheorem{Proposition}[Theorem]{Proposition}
\newtheorem{Remark}[Theorem]{Remark}

\def\m{\mathcal {M}}

\makeatletter \@addtoreset{equation}{section}
\title[Optimal transportation for generalized Lagrangian] {Optimal transportation for generalized Lagrangian}

\author[J.Li]{Ji Li}

\address{Department of mathematics, Nanjing University, 22 Hankou Road, 210093 Nanjing, China.}
\email{adailee.hepburn@gmail.com}

\author[J-L. Zhang]{Jianlu Zhang}
\address{Department of mathematics, Nanjing University, 22 Hankou Road, 210093 Nanjing, China.}
\email{mg0921036@smail.nju.edu.cn}

\begin{document}
\maketitle

\begin{abstract}
In this paper, we study the optimal transportation for generalized Lagrangian $L=L(x, u,t)$, and consider the cost function as following:
$$c(x, y)=\inf_{\substack{x(0)=x\\x(1)=y\\u\in\mathcal{U}}}\int_0^1L(x(s), u(x(s),s), s)ds.$$
Where $\mathcal{U}$ is a control set, and $x$ satisfies the following ordinary equation:
$$\dot{x}(s)=f(x(s),u(x(s),s)).$$
We prove that under the condition that the initial measure $\mu_0$ is absolutely continuous w.r.t. the Lebesgue measure, the Monge problem has a solution, and the optimal transport map just walks along the characteristic curves of the corresponding Hamilton-Jacobi equation:
\begin{equation*}
\begin{cases}
V_t(t, x)+\sup_{\substack{u\in\mathcal{U}}}<V_x(t, x), f(x, u(x(t), t),t)-L(x(t), u(x(t), t),t)>=0.\\
V(0,x)=\phi_0(x)
\end{cases}
\end{equation*}
\end{abstract}

\begin{quote}
\footnotesize {\it Keywords}: optimal control, Hamilton-Jacobi equation, characteristic curve, viscosity solution, optimal transportation,  Kantorovich pair, initial transport measure.\\

\end{quote}

\section{Introduction}
Given a pile of soil and an excavation that we want to fill up with the soil. As early as 1781, Monge posed the question to find an optimal way to do this work. This is the primary statement of Monge problem. We can model the pile of soil and the excavation by two probability measures as the density of the mass. Concretely, given a space $M$, and a continuous function, called the cost function $c(x,y):M\times M\rightarrow  \mathrm{R}$, given two probability measures $\mu_{0}$ and $\mu_{1}$ on $M$, find the mapping $t:M\rightarrow M$ which transport $\mu_{0}$ to $\mu_{1}$ and minimize the total cost $\int_{M}c(x,t(x))d\mu_{0}$. But unfortunately, the Monge problem is not always well-posed, for instance, consider $\mu_{0}=\delta_{x},\mu_{1}=\frac{1}{2}(\delta_y+\delta_z) (x, y, z\in M)$, the Monge problem has no solution since there is no map $t$ such that $t_{\sharp}\mu_{0}=\mu_{1}$.

It is difficult to study the Monge problem directly. In 1942, Kantorovich raised a generalized problem, later called the Monge-Kantorovich (MK) problem, we will state it concretely in section 4. As proved in \cite{K1}, (MK) problem always has a solution. It turned out that through the solution of (MK) problem, sometimes we can construct the solution of the Monge problem. Recently, several mathematicians studied the connections between Aubry-Mather-Fathi theory, optimal transportation problem and Hamilton-Jacobi equations \cite{DGG}\\\cite{EG}\cite{G}\cite{ P}\cite{V}\cite{ W}.

In \cite{BB}, they considered a Lagrangian function $L(x,v,t):TM\times[0, T]\rightarrow \mathrm{R}$ which satisfies the Tonelli conditions introduced by Mather, and considered the cost by:
   $$c(x,y)=\min_{\substack{\gamma}}\int_{0}^{T}L(\gamma(t),\dot{\gamma}(t),t)dt,$$
 where the minimum is taken on the set of absolutely continuous curve $\gamma:[0,1]\rightarrow M$ with $\gamma(0)=x$ and $\gamma(1)=y$. One of the main results of \cite{BB} is that using the dual problem (MK), they proved the optimal transformation can be performed by a Borel map with the assumption that the initial measure $\mu_{0}$ is absolutely continuous w.r.t. the Lebesgue measure, thus the Monge problem has a solution.

 This paper is stimulated by \cite{BB}, but we use another approach called `optimal control', studying the optimal transportation problem of a generalized Lagrangian $L=L(x, u(x, t), t)$. Consider the cost function as following:

 $$c(x,y)=\inf_{\substack{x(0)=x\\x(1)=y\\u\in \mathcal{U}}}\int_{0}^{1}L(x(t),u(x(t),t))dt,$$
 where $\mathcal{U}$ is called a control set, which we will state more clearly in section 3.  $x(t)$ satisfies the following equation:
 $$\dot{x}=f(x(t),u(x(t),t)),$$
\noindent Here $f$ satisfies some regular condition, we will state it clearly in section 3.

 We start from the point of optimal control to study the original optimal transportation problem, called the Monge problem:
 $$\hspace{-4.0cm}\mathbf{(M)}      \ \ \ \ \ \ \ \ C(\mu_{0},\mu_{1})=\min_{\substack{t}}\{\int_{M}c(x,t(x))d\mu_{0}\ | \; t_{\sharp}\mu_{0}=\mu_{1}\}.$$
Under the same assumption in \cite{BB} that the initial measure $\mu_{0}$ is absolute continuous w.r.t. the Lebesgue measure, we prove that the corresponding Monge problem has a Borel map as its solution.

\section{Method of characteristic}
In this section, we briefly introduce the method of characteristics  in constructing a local classical solutions of the Cauchy problem for Hamilton-Jacobi equation like
\begin{equation}\label{(2.1)}
\begin{cases}
  \frac{\partial u}{\partial t}+H(x, \nabla u(t, x))=0,\ \ \   (t, x)\in[0,\infty)\times M,\\
  u(0, x)=u_0(x),  \ \ \ x\in \mathrm{R}^n.
\end{cases}
\end{equation}
where $H$ and $u_0$ are of $C^2$.

We suppose a priori that we have a solution $u\in C^2([0,T]\times M)$ of the above equation. We call the solution of the following equation:
\begin{equation}
\begin{cases}
  \dot{X}=H_p(X,\nabla u(t,X)), \\  X(0)=z,
\end{cases}
\end{equation}
$t\rightarrow (t,X(t,z))$ a characteristic curve associated to $u$  starting from $z$.

Now we set
$$U(t,z)=u(t,X(t,z)),\ \ \  P(t,z)=\nabla u(t,X(t,z)).$$
Easy calculation shows that:
\begin{equation}
  \begin{cases}
    \dot{U}=-H(X,P)+P\cdot H_p(X,P),\\
    \dot{p}=\nabla u_t(t,x)+H_x(x,\nabla u(t,x))+\nabla^2u(t,x)H_p(x,\nabla u(t,x))
  \end{cases}
\end{equation}
 Since we have:
 $$0=\nabla (u_t(t,x)+H(x,\nabla u(t,x)))=\nabla u_t(t,x)+H_x(x,\nabla u(t,x))+\nabla^2u(t,x)H_p(x,\nabla u(t,x))$$
 we obtain that:
 $$\dot{P}=-H_x(X,\nabla u(t,X))=-H_x(X,P).$$
 Thus, the pair$(X,P)$ solves the following ordinary differential equation:
 \begin{equation}\label{(2.4)}
   \begin{cases}
     \dot{X}=H_p(X, P)\\
     \dot{P}=-H_x(X, P)
   \end{cases}
 \end{equation}
 with initial condition $X(0)=z, P(0)=\nabla u_0(z)$. While $U$ satisfies:
 \begin{equation}\label{(2.5)}
   \dot{U}=-H(X, P)+P\cdot H_p(X, P), \ \ \ U(0, z)=u_0(z).
 \end{equation}
Obviously, $X, P$ and $U$ are uniquely determined by the initial value $u_0$. The above arguments suggests that we can obtain a solution to the Hamilton-Jacobi equation (\ref {(2.1)}) by solving the characteristic system (\ref{(2.4)}) provided the map $z\rightarrow X(t, z)$ is invertible. In \cite{CS}, they list a classic result:

\begin{Theorem}\label{(*)}
For any $z\in \mathbb{R}^n$, let $X(t, z), P(t, z)$ denote the solution of the characteristic system (\ref{(2.4)}) and let $U(t, z)$ be defined by (\ref{(2.5)}). Then, there exists $T^*\geq 0$ such that the map $z\rightarrow X(t, z)$ is invertible with $C^1$ inverse $x\rightarrow Z(t, x)$, and there exists a unique solution $u\in C^2([0, T^*)\times M))$ of (\ref{(2.1)}), which is given by
\begin{equation}
u(t, x)=U(t; Z(t;x)), \ \ \  (t, x)\in [0, T^*)\times M.
\end{equation}
\end{Theorem}
\begin{Remark}
The $T^*$ in Theorem (\ref{(*)}) is the uniform maximum time such that any two characteristic curves don't intersect with each other when the time $T< T^*$.
\end{Remark}

\section{Optimal control problem for Bolza problem}

Optimality is a universal principle of life. The first basic ingredient of an optimal control problem is the so-called control system, it gives the possible behaviors. Usually, the control system is described by ordinary differential equations, in this paper, we consider the following form:
$$\quad\quad\quad\quad\quad\quad\quad\quad \dot{x}=f(x,u), \quad\quad\quad x(t_{0})=x_{0}. \;\;\;\;\;\;\;\;\;\;\;\quad\quad\quad\quad  \text{(ODE)}$$

\noindent here, $x$ is the state, $x\in M$, $t$ is the time, $x_{0}$ is the initial state, $t_{0}$ is the initial time, $u$ is called the control which depends on $x, t$, i.e. $u=u(x,t\in \mathbb{R})$.

Here are some basic assumptions on the control system:\\
(A1) There exists $K_0>0$ such that $|f(x,u)|\neq K_0(1+|x|+|u|), \forall x\in\mathcal{R}^n, u\in \mathcal{U}$;\\
(A2) $f$ is Lipschitz with respect to $x$: there exists $K_1>0$ such that $|f(x_2, u(x_2,t))-f(x_1, u(x_1,t))|\leq K_1|x_2-x_1|$, for all $x_1, x_2\in M, u\in \mathcal {U}$;\\
(A3) $f$ is $C^{1,1}$ with respect to $x$: there exists $K_2>0$ such that $||f_x(x_2, u(x_2,t)))-f_x(x_1, u(x_1,t))||\leq K_2|x_2-x_1|$, for all $x_1, x_2\in M, u\in \mathcal{U}$.

The assumption (A2) ensures that the existence of a unique global solution to the control system.\\

The second basic ingredient is the cost functional, it associates a cost with each possible behavior. An optimal control problem consists of choosing a control $u$ in the (ODE)in order to minimize the cost functional. 

Let $L=L(x, u(x,t),t)$ satisfies the following conditions:\\
(L1) $L$ is superlinear w.r.t $u$;\\
(L2) $L$ is locally Lipschitz w.r.t $x$: $\forall R>0$, there exists $\alpha_R$ such that:
$$|L(x_2,u(x_2,t),t)-L(x_1,u(x_1,t),t)|\leq \alpha_R|x_2-x_1|, \forall x_1, x_2 \in B_R;$$
$$|L(x, u_1(x,t),t)-L(x, u_2(x,t),t)|\leq C|u_1-u_2|, \forall x \in M;$$
(L3) $\forall x\in M$, the following set is convex:\\

\ \ $\mathcal{L}=\{(l,v)\in\mathcal{R}^{n+1}$: there exists $u\in \mathcal{U}$ such that $v=f(x,u), l\leq L(x,u)\}.$\\

For the convenience and unity, we consider the following optimal control problem:
  $$ \hspace{-4.5cm}(\mathbf{BP})\;  J(t,x,u)=\int_{0}^{t}L(x(s),u(x(s),s))ds+I(x(0)).$$
Where $L\in C^2(M\times\mathbb{R}\times[0, T]), I\in C(M)$, and $x$ is the solution of the control system:
  \begin{equation}\label{(3.1)}
   \begin{cases}
     \dot{x}=f(x(s),u(x(s),s)),\\
     x(t)=x
   \end{cases}
 \end{equation}
 Note that fixing endpoint is the same as fixing initial point in essential. Here we use the same symbol $``x"$ to denote a point or a solution of (\ref{(3.1)}) without confusion.
 Now, we consider a time period of [0, T], and let $t$ range over [0,T],  $x$ range over $M$. We introduce the value function:
\begin{Theorem}\cite{CS}
Assume (A1), (A2), (L1), (L2), (L3) hold, and the initial cost $I$ is continuous. Then,$\forall  (x, t)\in M\times[0,T]$, there exists an optimal control $u\in\mathcal{U}$  for $\mathbb(BP)$.

\end{Theorem}

 \begin{Definition}
   Given $(t,x)\in [0,T]\times M$, define:
   $$V(t,x)=\inf\{ J(t,x,u)\; |\; u:[0,t]\times M\rightarrow \mathcal{U}\},$$
   the function $V$ is called the value function of the control problem ($\mathbf{BP}$).
\end{Definition}

  One of the most important principle in optimal control problem is the so called `` principle of optimality" for $V(t,x):$
   \begin{Theorem}
     (Principle Of Optimality): For every $(t,x)\in [0,T]\times M$ and $\Delta t\in(0,t)$, the value function $V$ satisfies the following relation:
     \begin{equation}\label{(3.2)}
     V(t,x)=\inf_{\substack{u_[t-\Delta t,t]}}\{\int_{t-\Delta t}^{t}L(x(s),u(x(s),s))ds+V(t-\Delta t,x(t-\Delta t))\}.
     \end{equation}
   \end{Theorem}

 The above principle of optimality means that in order to search for an optimal control, we can search over a small time interval for a control that minimizes the cost over this interval plus the subsequent optimal cost-to-go.

 Now, let's look at the value function $V(t,x)$ in a deep sight. We assume $V(t,x)\in C^1([0,T]\times M)$, consider the right-hand side of (\ref{(3.2)}):
 \begin{equation}\label{(3.3)}
   x(t-\Delta t)=x-f(x,u)\Delta t+o(\Delta t),
 \end{equation}
so, we can express $V(t-\Delta t,x(t-\Delta t))$ as follow:
\begin{equation}
\begin{split}
  &V(t-\Delta t,x(t-\Delta t))\\
  =&V(t,x)+V_t(t,x)(-\Delta t)+<V_x(t,x),f(x,u)(-\Delta t)>+o(\Delta t),
  \end{split}
\end{equation}
 and we have:
\begin{equation}\label{(3.5)}
  \int_{t-\Delta t}^{t}L(x(s),u(x(s),s))ds=L(x(t),u(x(t),t))\Delta t+o(\Delta t).
\end{equation}
 Substituting (\ref{(3.3)}) and (\ref{(3.5)}) into the right-hand side of (\ref{(3.2)}), we get:
 \begin{equation}
 \begin{split}
   V(t,x)=\inf_{\substack{u_{[t-\Delta t, t]}}}\{&L(x(t),u(x(t),t))\Delta t+V(t,x)+V_t(t,x)(-\Delta t)\\
   &+<V_x(t,x),f(x,u)(-\Delta t)>+o(\Delta t)\}.
   \end{split}
 \end{equation}
   After simple calculation, we obtain:
   \begin{equation}
     0=\inf_{\substack{u_{[t-\Delta t, t]}}}\{ L(x(t),u(x(t),t))-V_t(t,x)
    -<V_x(t,x),f(x,u)>\}.
   \end{equation}
 This is equivalent to the following equation:
 \begin{equation}\label{(3.8)}
   V_t(t,x)+\sup\{<V_x(t,x),f(x,u)>-L(x(t),u(x(t),t))\}=0.
 \end{equation}
It is remarkable to note that there are several forms which are equivalent to (\ref{(3.8)}), but (\ref{(3.8)}) is the reasonable form since flipping the sign in a PDE affects its viscosity solutions, it will turn out that the above sign convention is the correct one.
Let
\begin{equation}\label{(3.9)}
H=\sup_{\substack{u\in\mathcal{U}}}<p, f(x,u)>-L(x,u,t).
\end{equation}
We have the following theorom:
\begin{Theorem}\label{3.4}
If the initial condition is equal to the initial cost, i.e. $u_0(x)=I(x)$, then V(t, x) is the unique viscosity solution of (\ref{(2.1)} ) with Hamiltonian function $H$ as in (\ref{(3.9)}).
\end{Theorem}
Proof:
Step 1.(viscosity subsolution) For every $C^1$ test function $\varphi=\varphi(x,t)$, assume $\varphi\geq V$, and $\varphi-V$ attains a local minimum at $(t_0,x_0)$, we need to prove:
$$\varphi(t,x)+\sup_{\substack{u\in \mathcal{ U}}}\{<\varphi_x(t,x),f(x,u)>-L\}\leq 0.$$
Suppose the contrary, there exists a $C^1$ function $\varphi$, and a control $u_0$, s.t:
$$\varphi(t_0,x_0)=V(t_0,x_0), \varphi(t,x)=V(t,x) \forall (t,x)\in U_\delta(t_0,x_0).$$
and
$$\varphi_t(t-0,x_0)+<\varphi_x(t_0,x_0),f(x,u)>-L>0.$$
Since
\begin{equation}
\begin{split}
&V(t_0-\Delta t,x(t_0-\Delta t))-V(t_0,x_0)\\
  \leq&\varphi(t_0-\Delta t,x(t_0-\Delta t))-\varphi(t_0,x_0)\\
=&-\int_{t_0-\Delta t}^{t_0}\frac{d\varphi(t,x)}{dt}\\
=&-\int_{t_0-\Delta t}^{t_0}\varphi_t(t,x)+<\varphi_x(t,x),f(x,u)>\\
\leq&-\int_{t_0-\Delta t}^{t_0}Ldt.
\end{split}
\end{equation}
This is contradictory to the principle of optimal, so $V(t,x)$ is a viscosity subsolution.\\
 Step 2.(viscosity supersolution) For every $C^1$ test function $\varphi=\varphi(t,x)$, assume $\varphi\leq V$,and $\varphi-V$ attains a local minimal at $(t_0,x_0)$, we need to prove:
$$\varphi(t,x)+\sup_{\substack{u\in \mathcal{ U}}}\{<\varphi_x(t,x),f(x,u)>-L\}\geq 0.$$
Suppose the contrary, there exists a $C^1$ function $\psi$, and a control $u'_0$,s.t:
$$\psi(t_0,x_0)=V(t_0,x_0),\psi(t,x)\leq V(t,x),\forall (t,x)\in U_\delta (t_0,x_0)$$
and
$$\psi_t+<\psi_x,f>-L<0.$$
Since
\begin{equation}
\begin{split}
\int_{t_0-\Delta t}^{t_0}L=&V(t_0,x_0)-V(t_0-\Delta t,x(t_0-\Delta t))\\
  \leq &\psi(t_0,x_0)-\psi(t_0-\Delta t,x(t_0-\Delta t))\\
  =&\int_{t_0-\Delta t}^{t_0}\frac{d\psi(t,x)}{dt}=\int_{t_0-\Delta t}^{t_0}\psi_t(t,x)+<\psi_x,f(x,u)>dt\\
  <&\int_{t_0-\Delta t}^{t_0}L(x,u)dt.
\end{split}
\end{equation}
We get a contradictory, thus $V(t,x)$ is a viscosity supersolution.\\
Step 3.(Uniqueness) By the comparison principle, the viscosity solution is unique!

From the above three steps, we have proved that $V(t,x)$ is the unique viscosity solution.

Let the Hamiltonian $H$ be as (\ref{(3.9)}), we have the following theorem:
\begin{Theorem} \cite{CS}
 Let $(u,x)$ be an optimal pair for the point $(t_0,x_0)\in[0, T]\times M$ and let $p:[0, t]\rightarrow \mathrm{R}^n$ be a dual arc associated with $(u,x)$, then $(x,p)$ solves the system:
  \begin{equation}
  \begin{cases}
    \dot{x}(s)=H_x(x(s),p(s)),\\
    \dot{p}(s)=-H_p(x(s),p(s))
  \end{cases}
  \end{equation}
  for all $s\in[0, t]$. As a consequence, $x, p$ are of class $C^1$.
\end{Theorem}

\begin{Remark}
 The above theorem tells us that the optimal trajectory of the optimal control problem is nothing else, but just the characteristic curve of the corresponding Hamilton-Jacobi equation.
\end{Remark}

\section{Optimal transportation related to $L(x,u,t)$}
\noindent $\mathbf{Statement\  of \ the\  Monge\  problem:}$ Let $M$ be a compact manifold, and a continuous cost function $c(x,y)$ is given, $c(x,y):M\times M\rightarrow \mathrm{R}$, $\mu_0,\mu_1 $ are two probability measures on $M$. Find the mapping $\Psi:M\rightarrow M$ which transport $\mu_0$ into $\mu_1$, and minimize the total cost: $$\int_Mc(x,\Psi(x))d\mu_0,$$
where the ``transport" means ``push-forward", which is defined as following :

\begin{Definition}
(Push-forward) Let $t:X\rightarrow Y$ be a measurable map, $\mu$ is a measure, define:
$$t_\sharp\mu(B)=\mu(t^{-1}(B))$$
for every Borel set $B\subset Y.$
\end{Definition}

The Monge problem can be stated as:
$$\hspace{-2.4cm}\mathbf{(M)}:=\min\{\int_{X\times Y}c(x,\Psi(x))d\mu_0|\Psi:M\rightarrow M, \Psi_\sharp\mu_0=\mu_1\}.\ \ \ \ \ \ \ \ \ \ $$(
In the rest of this paper, we consider the cost function as following:
\begin{equation}
  c(x,y)=\inf_{\substack{x(0)=x\\x(1)=y\\u\in\mathcal{U}}}\int_{0}^{1}L(x(s),u(x(s),s))ds.
\end{equation}
where $x$ satisfies the following ordinary equation:
\begin{equation}
\dot{x}(s)=f(x(s),u(x(s),s)).
\end{equation}
As we said in the introduction, the Monge problem is not always well-posed, since sometimes there is no map $\Psi$ such that $\Psi_\sharp\mu_0=\mu_1$. In the study of the Monge problem, mathematicians find that it is a good and effective approach to consider the (MK) problem first:\\
\noindent $\mathbf{Statement\ of\ the\ Monge-Kantorovich\ problem:}$

$$( \mathbf{MK} ):= \min\{\int_{X\times Y}c(x,y)d\gamma(x,y)|\gamma\in \mathcal{P}(X\times Y), (\pi_1)_\sharp\gamma=\mu_0, (\pi_2)_\sharp\gamma=\mu_1\}.$$

Under some proper conditions, we can use the solution of the (MK) problem to construct the solution of the Monge problem. In this paper, we use this approach.

The measure $\gamma$ which satisfies $(\pi_1)_\sharp\gamma=\mu_0, (\pi_2)_\sharp\gamma=\mu_1$ in ($\mathbf{MK}$) is called a transport plan. We denote the transport plans by $\mathcal{K}(\mu_0,\mu_1)$. Observe that if $\Psi$ is admissible for Monge problem, i.e. $\Psi_\sharp\mu_0=\mu_1$, then the measure $\gamma=(Id\times\Psi)_\sharp \mu_0$ is a transport plan for ($\mathbf{MK}$).Moreover, the class of transport plans is always non-empty since it always contains $\mu_0\times\mu_1$. As proved in \cite{K1}, the ``min" in ($\mathbf{MK}$) always can be achieved by a transport plan, and we call the transport plan which realizes the ``min" an optimal transfer plan \cite{BB}.

Let the total cost
\begin{equation}
  C(\mu_0,\mu_1)=\min_{\substack{\gamma\in\mathcal{K}(\mu_0,\mu_1)}}\int c(x, y)d\gamma.
\end{equation}
\begin{Definition}
  (Admissible Kantorovich Pair) A pair of continuous function $(\phi_0,\phi_1)$ is called an admissible Kantorovich pair if it satisfies the following relations:
  \begin{equation}
  \phi_1(x)=\min_{\substack{y\in M}}\phi_0(y)+c(x, y),\  and \ \phi_0(x)=\max_{\substack{y\in M}}\phi_1(y)-c(x, y),
  \end{equation}
  for all point $x \in M$.
\end{Definition}

Kantorovich proved the following notable result:
\begin{Theorem}
\begin{equation}
 C(\mu_0,\mu_1)=\max_{\substack{(\phi_0,\phi_1)}}(\int_M\phi_1d\mu_1-\int_M\phi_0d\mu_0). \end{equation}
 where the ``max" is taken on the set of admissible Kantorovich pair $(\phi_0,\phi_1)$.
\end{Theorem}
This is the so-called dual Kantorovich problem. The admissible pairs which attain the ``max" are called optimal Kantorovich pairs, usually, it is not unique!
\begin{Proposition}
  If $\gamma$ is an optimal transfer plan , and $(\phi_0,\phi_1)$ is an optimal Kantorovich pair, then the support of $\gamma$ is contained in the following set:

  $$\{(x, y)\in M\times M such\  that\
    \phi_1(y)-\phi_0(x)=c(x, y)\}\subset M\times M$$
\end{Proposition}
\begin{Proposition}
  We have:
  $$c(x, y)=\max_{\substack{(\phi_0, \phi_1)}}\phi_1(y)-\phi_0(x).$$
  where the ``max" is taken on the set of admissible Kantorovich pairs.
\end{Proposition}
It is deserved to note that for different pair of $(x, y)$, the ``max" may be achieved by different admissible Kantorovich pair.\\

Let $(\phi_0, \phi_1)$ be an optimal Kantorovich pair,  and
$$c_0^t(x, y)=\inf_{\substack{x(0)=x\\x(1)=y\\u\in\mathcal{U}}}\int_{0}^{t}L(x(s),u(x(s),s))ds.$$
we construct a function on $M\times [0, T]$ as following:
\begin{equation}
  U(x,t)=\min_{\substack{y\in M}}\phi_0(x)+c_0^t(x, y),
  \end{equation}

Recall the definition of the value function for $\mathbf{(BP)}$:
$$V(t, x)=\inf_{\substack{u\in\mathcal{U}}}\int_0^tL(x(s),u(x(s),s),s)ds+I(x(0)),$$
  where $x$ satisfies the control system (\ref{(3.1)}). Obviously,  we can see easily that $U(x,t)$ is just the same as the value function $V(t,x)$ with $I(x)=\phi_0(x).$

  So, the total cost $C(\mu_0,\mu_1)$ can be realized by a viscosity solution of the Hamilton-Jacobi equation (\ref{(2.1)}) with a proper initial condition. We have:
  \begin{equation}
    C(\mu_0,\mu_1)=\int_M V(1, x)d\mu_1-\int_M V(0, x)d\mu_0.
  \end{equation}

  Let's stop here for a little while, and look at the value function from another point of view. Assume $u^*$ realizes the ``inf" in the definition of value function $V(t,x)$, we define an operator $\mathcal{T}: C^{ac}([0,t],R)\rightarrow C^{ac}([0,t],R)$  as following:

\begin{equation}
\mathcal{T}_t \phi_0(x)=\phi_0(\gamma(0))+\int_0^tL(\gamma(s),u^*(\gamma(s),s),s)ds.
\end{equation}

Obviously, $\mathcal{T}_t \phi_0(x)=V(t,x)$.

As in \cite{SY,WY}, we call $\mathcal{T}_t$ the solution semigroup since it is determined by the viscosity solution of the Hamilton-Jacobi equation (\ref{(2.1)}). Actually, $\mathcal{T}_t$ is indeed a semigroup.
\begin{Theorem}(semi-group property)
$\{\mathcal{T}_t\}_{t\geq 0}$ is a one-parameter semigroup from $C(M, \mathbf{R})$.
\end{Theorem}
Proof: $\mathcal{T}_s\phi_0(x)=u(x,s)$, this means that the operator $\mathcal{T}_s$ just sends the initial value $\phi_0$ into the corresponding viscosity solution at time $s$ under this initial condition. Similarly, $\mathcal{T}_{t+s}\phi_0(x)=u(x,s+t)$ means that the operator  $\mathcal{T}_{t+s}$ sends the initial value $\phi_0$ into the corresponding viscosity solution at time $t+s$. Attention $\mathcal{T}_t\circ \mathcal{T}_s\phi_0(x)$ means that the operator $\mathcal{T}_t$ sends the initial value $\mathcal{T}_s\phi_0$ to the corresponding viscosity solution at time $t$. Since the viscosity solution of Hamilton-Jacobi (\ref{(2.1)}) is unique, it is obvious that $\mathcal{T}_t\circ\mathcal{T}_s\phi_0=\mathcal{T}_{t+s}\phi_0$, for arbitrary initial value $\phi_0\in C(M, \mathbf{R})$.

\begin{Definition}
  Let $(\phi_0,\phi_1)$ denote an optimal Kantorovich pair, and $u^*$ is one of the control which reaches the ``inf" in the definition of the value function $V$. We define $\mathcal{F}(\phi_0,\phi_1) \subset C^2([0, 1], M)$ the set of curves $\gamma(t)$ such that:
  \begin{equation}\label{(4.9)}
    \phi_1(\gamma(1))=\phi_0(\gamma(0))+\int_0^1 L(\gamma(t),u^*(\gamma(t),t))dt.
  \end{equation}
\end{Definition}
\begin{Remark}
  Obviously,  if we rewrite (\ref{(4.9)}) into the following form:
 \begin{equation}\label{4.10}
V(1, \gamma(1))=V(0,\gamma(0))+\int_0
^1 L(\gamma(t), u^*(t, \gamma(t), t)dt.
\end{equation}
We can see that $\mathcal{F}(\phi_0,\phi_1) $ is just the set of pieces of characteristic curves.
\end{Remark}

Let $\mathcal{F}_0(\psi_0, \psi_1)$ be the set of initial state $(x, p)\in T^*M$ such that the curve $t\rightarrow \pi \circ\psi_0^t(x, v)$ belongs to $\mathcal{F}(\phi_0, \phi_1)$. We have the following lemma:
\begin{Lemma}
  The maps $\pi$ and $\pi\circ \psi_0^1:\mathcal{F}_0(\phi_0, \phi_1)\rightarrow M$ are surjective. If $x$ is a differentiable point of $\phi_0$, then the set $\pi^{-1}(x)\cap\mathcal{F}(\phi_0,\phi_1)$ contains only one point. There exists a Borel measurable set $\Sigma\subset M$ of full measure, whose points are differentiable points of $\phi_0$, and such that the map:
  $$x\rightarrow S(x)=\pi^{-1}(x)\cap\mathcal{F}(\phi_0,\phi)$$
  is Borel measurable on $\Sigma$.
\end{Lemma}
Proof: For each $x\in M$, there exists a characteristic curve such that (\ref{(4.9)}) and (\ref{4.10}) are satisfied, so the projection $\pi\circ\psi_0^1$ from $\mathcal{F}_0(\phi_0, \phi_1)$ to $M$ is surjective. For $\pi$, it's similar! 

Next we consider a differentiable point $x$ of $\phi_0$. By the characteristic method introdeced in section 2,  the characteristic curves don't cumulate together at $x$, i.e. there is only one characteristic curve start from $x$. We construct $S$ as following:
\begin{equation}
S(x)=\pi^{-1}(x)\cap\mathcal{F}(\phi_0,\phi),
\end{equation}
Since $\phi_0$ is Lipschitz, the set of differentiable points of $\phi_0$ is of total Lebesgue measure. Notice that there exists a sequence of compact sets $K_n$ such that $\phi_0$ is differentiable at each point of $K_n$, and the Lebesgue measure of $M-K_n$ is converging to zero. For each $n$, the set $\pi^{-1}(K_n)\cap \mathcal{F}(\phi_0,\phi_1)$ is compact, and the canonocal projection $\pi$ restricted to this set is injective and continuous, so the inverse function $S$ is continuous on $K_n$. And as a consequence, $S$ is a Borel measurable map on $\Sigma=\cup_n K_n$. \\

Now let
$$H(x, p)=\max_{\substack{u}}{<p, f(x, u)>-L(x, u, t)}.$$
And let $m_0\in \mathcal{B}(T^*M)$ be a Borel probability measure on the cotangent bundle $T^*M$. We call $m_0$ an initial transport measure if the measure $\eta$ on $M\times M$ given by
$$\eta=(\pi\times(\pi\circ \psi))_\sharp m_0$$ 
is a transport plan. Here $\psi$ is the time one map of the Hamiltonian flow, and $\pi: T^*M\rightarrow M$ is the canonical projection. We denote $\mathcal{I}(\mu_0, \mu_1)$ the set of initial transport measures. Obviously, we can define the action of an initial transport measure as following:
\begin{equation}
A(m_0)=\int_{M\times M}c(x, y)d\eta.
\end{equation}

\begin{Lemma}
The mapping $ (\pi\times(\pi\circ\psi))_\sharp:\mathcal{I}(\mu_0, \mu_1)\rightarrow \mathcal{K}(\mu_0, \mu_1)$ is surjective.
\end{Lemma}
Proof: We shall prove that for arbitrary probability measure $\eta\in\mathcal{B}(M\times M)$, there exists a probability measure $m_0\in \mathcal{B}(T^*M)$ such that $(\pi\times(\pi\circ\psi_0^T))_{\sharp}m_0=\eta$. Since the set of probability measures on $M\times M$ is the compact convex closure of the set of Dirac probability measures, so we just need to prove it when $\eta$ is a Dirac probability measure. Assume $\eta$ be a Dirac probability measure supported at $(x,y)\in M\times M$, and $\gamma: [0,T]\rightarrow M$ be a curve with boundary conditions $\gamma(0)=x,\gamma(1)=y$ satisfying the control system. Let $m_0$ be a Dirac probability measure supported at $\gamma(0), \dot{\gamma}(0)$. Obviously, we have:
$$\eta=(\pi\times(\pi\circ\psi_0^T))_{\sharp}m_0$$.

\noindent Where $\eta$ is a transport plan which is determined by $m_0$ through $ (\pi\times(\pi\circ\psi))_\sharp$.
\begin{Theorem}
  If $\mu_0$ is absolutely continuous with respect to Lebesgue measure, then, there exists a Borel section $S: M\rightarrow TM$ such that the map $\pi\circ\psi_0^1\circ S$ is an optimal transport map between $\mu_0$ and $\mu_1$. The section $S$ is unique in the sense of $\mu_0$ -almost everywhere.
\end{Theorem}

Proof: Since $\Sigma$ is of full Labesgue measure, consider $m_0=S_\sharp(\mu_0|_\Sigma)$, which is a probability measure on $T^*M$, concentrated on $\mathcal{F}(\phi_0, \phi_1)$, and $\pi_\sharp m_0=\mu_0$. Since $\pi$ induces an isomorphism between $\pi^{-1}(\Sigma)\cap \mathcal{F}_0(\phi_0, \phi_1)$, and by theorem ??, we have $C_0^1(\mu_0,\mu_T)=\min_{\mathcal{I}}A(m_0)$, and $m_0$ is the unique initial transport measure. Thus, $\m_0$ is the unique optimal transport measure. $\pi\circ\psi\circ S$ is an optimal transport map between $\mu_0$ and $\mu_1$.

\end{document}